\documentclass[12pt]{article}
\usepackage{amsthm,amsfonts,amssymb,amscd}

\begin{document}

\begin{center} \large \bf Birational rigidity of quartic three-folds
\\ with a double point of rank 3
\end{center}\vspace{0.5cm}

\centerline{A.V.Pukhlikov}\vspace{0.5cm}

\parshape=1
3cm 11cm \noindent {\small \quad\quad\quad \quad\quad\quad\quad
\quad\quad\quad {\bf }\newline We prove that a general
three-dimensional quartic $V$ in the complex projective space
${\mathbb P}^4$, the only singularity of which is a double point
of rank 3, is a birationally rigid variety. Its group of
birational self-maps is, up to the finite subgroup of biregular
automorphisms, a free product of 25 cyclic groups of order 2. It
follows that the complement to the set of birationally rigid
factorial quartics with terminal singularities is of codimension
at least 3 in the natural parameter space.

Bibliography: 15 items.} \vspace{1cm}


\vspace{1cm}

\section*{Introduction}

{\bf 0.1. The main result.} Let $V\subset {\mathbb P}^4$ be a
hypersurface of degree 4 in the complex projective space, which in
some affine chart ${\mathbb A}^4\subset {\mathbb P}^4$ with
coordinates $(z_1,z_2,z_3,z_4)$ is given by an equation
$f(z_1,\dots, z_4)=0$, where
$$
f=q_2(z_*)+q_3(z_*)+q_4(z_*),
$$
where $q_2(z_*)=z_1^2+z_2^2+z_3^2$ is a quadratic form of rank 3,
the condition
$$
q_3(0,0,0,1)\neq 0
$$
is satisfied and the set of common zeros of the homogeneous
polynomials $q_2$, $q_3$, $q_4$ on ${\mathbb
P}^3_{(z_1:z_2:z_3:z_4)}$ consists of 24 distinct points. Assume
that the point $o=(0,0,0,0)$ is the unique singularity of the
projective quartic $V$. It is easy to check that this singularity
is resolved by the blowing up of the point $o$, which implies at
once (see \S 1), that $V$ is a factorial variety with the only
terminal singular point and $\mathop{\rm Pic} V={\mathbb Z} H$,
where $H$ is the class of a hyperplane section, so that $V$ is a
Fano three-fold of index 1 (since $K_V=-H$). The aim of the
present paper is to prove the following claim.

{\bf Theorem 0.1.} {\it The variety $V$ is birationally rigid: if
$\chi\colon V\dashrightarrow V'$ is a birational map onto the
total space $V'$ of some Mori fibre space $\pi'\colon V'\to S'$,
then $S'$ is a point and $V'$ is isomorphic to $V$ in the
biregular sense.}

This claim has been stated in \cite{CoMe} as a conjecture and up
to this day remained an open problem. Let us emphasize that if
$V'$ is a three-dimensional Fano variety with ${\mathbb
Q}$-factorial terminal singularities and the Picard number 1,
which is birational to the quartic $V$, then Theorem 0.1 claims
that $V'\cong V$, but it {\it does not} claim that every
birational map between these varieties is a biregular isomorphism
(the quartic $V$ is not birationally superrigid): identifying $V'$
and $V$, we can interpret a birational map $\chi\colon
V\dashrightarrow V'$ as a birational self-map $\chi\in\mathop{\rm
Bir} V$. The group $\mathop{\rm Bir} V$ is very large, its
description is identical to the description of the group of
birational self-maps of the three-dimensional quartic with a {\it
non-degenerate} quadratic singularity \cite{Pukh89c}.

Given the assumptions about the polynomial $f(z_*)$, there are 24
lines $L_1,\dots, L_{24}\subset V$ passing through the singular
point $o$ (they correspond to the 24 common zeros of the
polynomials $q_2$, $q_3$, $q_4$). With each of these lines $L_i$ a
birational involution $\tau_i\in \mathop{\rm Bir} V$ is naturally
associated, $\tau_i^2=\mathop{\rm id}_V$, see \cite{Pukh89c}, this
construction is recalled in Subsection 1.3.

Let $\tau_0\in \mathop{\rm Bir} V$ be the Galois birational
involution, corresponding to the projection of the quartic $V$
from the point $o$ onto ${\mathbb P}^3$ (this projection is a
rational map of degree 2). Set
$$
\mathop{\rm B}(V)=\langle \tau_0, \tau_1, \dots, \tau_{24}\rangle
\subset \mathop{\rm Bir} V
$$
to be the subgroup, generated by these 25 involutions. Theorem 0.1
is complemented by the following claim.

{\bf Theorem 0.2.} (i) {\it The subgroup $\mathop{\rm B}(V)$ is
freely generated by the involutions $\tau_0$, $\tau_1$, \dots,
$\tau_{24}$,
$$
\mathop{\rm B}
(V)=\mathop{*}\limits_{i=0}^{24}\langle\tau_i\rangle,
$$
and is a normal subgroup of the group $\mathop{\rm Bir} V$.

{\rm (ii)} The group $\mathop{\rm Bir} V$ is an extension of the
finite subgroup of biregular automorphisms $\mathop{\rm Aut} V$ by
the normal subgroup} $\mathop{\rm B}(V)$:
$$
1\to \mathop{\rm B}(V)\to \mathop{\rm Bir} V\to \mathop{\rm Aut}
V\to 1.
$$

For a general polynomial $f(z_*)$ the group $\mathop{\rm Aut} V$
is trivial, so that $\mathop{\rm Bir} V$ is a free product of 25
cyclic groups of order 2.

Theorems 0.1 and 0.2 give a complete description of the birational
type of the quartic $V$ (which is identical to the description of
the birational type of a quartic with a non-degenerate double
point of general position \cite{Pukh89c}). The claim of Theorem
0.1 related to Mori fibre spaces $V'/S'$ over a base of positive
dimension can be simpler and clearer stated as follows: the
quartic $V$ can not be fibred by a rational map into rational
curves or rational surfaces. Using the traditional language, on
the variety $V$ there are no structures of a conic bundle and no
pencils, a general member of which is a rational
surface.\vspace{0.3cm}


{\bf 0.2. Singularities of linear systems.} Now let us state the
main technical claim which implies Theorems 0.1 and 0.2 by easy
standard arguments. Let $\Sigma$ be a mobile linear system on $V$
(that is, $\Sigma$ has no fixed components). By the factoriality
of the variety $V$ and the lefschetz theorem the system $\Sigma$
is cut out on $V$ by hypersurfaces of degree $n[\Sigma]\geqslant
1$ in ${\mathbb P}^4$, that is, $\Sigma\subset |n[\Sigma]\, H|$.

{\bf Theorem 0.3.} {\it Assume that for general divisor
$D\in\Sigma$ the pair
\begin{equation}\label{13.09.24.1}
\left(V,\frac{1}{n[\Sigma]}D\right)
\end{equation}
is not canonical. Then in the set of birational involutions
$\{\tau_0,\tau_1,\dots,\tau_{24}\}$ there is precisely one
$\tau_a$, $0\leqslant a\leqslant 24$, such that the strict
transform $(\tau_a)_*\Sigma$ of the linear system $\Sigma$ with
respect to $\tau_a$ satisfies the inequality
$$
n[(\tau_a)_*\Sigma]<n[\Sigma],
$$
whereas for any other involution $\tau_i$, $i\neq a$, the opposite
inequality holds:}
$$
n[(\tau_i)_*\Sigma]>n[\Sigma].
$$

In the statement of the theorem one should take into account that
$\tau_a=\tau_a^{-1}$. The procedure of changing $\Sigma$ by
$(\tau_a)_*\Sigma$ is traditionally called {\it untwisting}.
Obviously, the untwisting process can not be infinite, and this
gives both claims, of Theorems 0.1 and 0.2, at once, see \S 1.

The pair (\ref{13.09.24.1}) being non-canonical means that for
some divisor $E$ over $V$ (that is, for some divisor $E$ on a
non-singular variety $\widetilde{V}$, where $\varphi\colon
\widetilde{V}\to V$ is a birational morphism), losing dimension on
$V$ (that is, $\dim \varphi(E)\leqslant 1$), the Noether-Fano
inequality
\begin{equation}\label{21.09.24.1}
\mathop{\rm ord}\nolimits_E \Sigma>n[\Sigma]\cdot a(E)
\end{equation}
holds, where $a(E)=a(E,V)\geqslant 1$ is the discrepancy of $E$
with respect to $V$ and in the sense of explanations given above,
$\mathop{\rm ord}\nolimits_E \Sigma=\mathop{\rm ord}\nolimits_E
\varphi^*\Sigma$. The divisor $E$ over $V$, or, more precisely,
the discrete valuation $\mathop{\rm ord}\nolimits_E (\cdot)$ of
the field of rational functions ${\mathbb C}(V)$, is traditionally
called a maximal singularity of the system $\Sigma$.

For the centre of $E$ on $V$ we have the following options:

\noindent (0) it is the singular point $o\in V$,

\noindent (1) it is one of the 24 lines $L_1$, \dots, $L_{24}$,

\noindent (2) it is a non-singular point of the quartic $V$,

\noindent (3) it is a curve which is not one of the 24 lines.

In order to prove Theorem 0.3, we will show that the options (2)
and (3) do not realize, in the case (0) the exceptional divisor of
the blow up of the point $o$ is already a maximal singularity, and
a certain claim about the uniqueness of the centre of a maximal
singularity is true, the precise statement of which we put off
until \S 1. The exclusion of the option (2) goes back to the
classical paper \cite{IM}. A part of the work described above, a
plan of the remaining work and some discussion make the contents
of \S 1.\vspace{0.3cm}


{\bf 0.3. The structure of the paper.} The paper is organized in
the following way. In \S 1 we give the construction of the 25
birational involutions $\tau_i$ and show that Theorems 0.1 and 0.2
can be easily derived from Theorem 0.3, and Theorem 0.3 itself
follows from some fact about the multiplicities of the linear
system at the point $o$ and along the 24 lines passing through
this point (Theorem 1.1). \S 2 contains some preparatory work for
the proof  of Theorem 1.1 and the statement of the main fact
(Theorem 2.1), which easily implies Theorem 1.1. Finally, \S 3
contains the proof of Theorem 1.1, the meaning of which is that if
the pair (\ref{13.09.24.1}) is not canonical, canonical outside
the point $o$, then the exceptional quadric $Q$ of the blow up of
the singular point is already a maximal singularity of the system
$\Sigma$. This fact is of key importance.\vspace{0.3cm}


{\bf 0.4. The birational geography of three-dimensional quartics.}
The birational superrigidity of smooth three-dimensional quartics
was shown in the classical paper \cite{IM}. The natural parameter
space for three-dimensional quartics is the 70-dimensional space
${\cal P}_{4,5}$ of homogeneous polynomials of degree 4 in 5
variables (more precisely, its projectivization), from which one
should remove an irreducible closed subset of dimension 38 (a
general polynomial in it is a product of a linear form and a cubic
form) and an irreducible closed subset of dimension 29 (a general
polynomial in which is a product of two quadratic forms). The
theorem of Iskovskikh and Manin gives a complete description of
the birational type of any quartic outside the (irreducible)
hypersurface of singular quartics.

The paper \cite{Pukh89c} proved the birational rigidity and gave a
complete description of the group of birational self-maps for
general three-dimensional quartics with a unique double point of
the maximal rank 4. The hardest part of that work (\cite[\S
5]{Pukh89c}) was to prove an analog of Theorem 1.1, see the end of
Subsection 0.3 above: if a linear system has a maximal singularity
over the quadratic point, then already the point itself ia maximal
for that system. Corti noticed \cite{Co00}, that this claim by the
inversion of adjunction follows immediately from the fact that a
pair $({\mathbb P}^1\times {\mathbb P}^1,\frac{1}{m} C)$ is log
canonical, where $C$ is a curve of bi-degree $(m,m)$ on ${\mathbb
P}^1\times {\mathbb P}^1$ (the latter is not difficult to show;
see, for instance, \cite[\S 7]{Pukh89c}; there are other proofs).
This simplification made it possible to prove the birational
rigidity \cite{Me04} and describe the group of birational
self-maps \cite{Me04,Shr08b} for any three-dimensional quartic
with non-degenerate double points in the assumption that is
${\mathbb Q}$-factorial. If we take into account that in
\cite{Ch2010} it was shown that a quartic with $\leqslant 8$
non-degenerate double points is factorial, and in
\cite{Shramov2007} it was shown that a quartic with $\leqslant 12$
non-degenerate double points (except for the case when there are
precisely 12 singular points and the quartic contains a
two-dimensional quadric) is ${\mathbb Q}$-factorial, we see that
the complete description of the birational type of
three-dimensional quartics with non-degenerate double points has
progressed rather far, up to codimension 11 in the parameter space
${\mathbb P}({\cal P}_{4,5})$.

However, a three-dimensional quartic can have worse singularities
than non-degenerate double points. The paper \cite{CoMe}
successfully studied the birational geometry of general quartics,
the affine equation of which in the notations of Subsection 0.1
satisfies the equality $\mathop{\rm rk} q_2=2$: it was shown that
such quartic has no structures of a conic bundle or fibration into
rational surfaces, and if it is birational to a Fano three-fold
$Y$ with ${\mathbb Q}$-factorial terminal singularities, then for
$Y$ there are only two options: either $Y$ is the quartic itself,
or $Y$ is a quasi-smooth complete intersection of type $3\cdot 4$
in the weighted projective space ${\mathbb P}(1^4,2^2)$, the
equation of which is determined by the equation of the original
quartic. In other words, in each of the two cases $Y$ is
determined up to an isomorphism. The proof makes use of the
Sarkisov program \cite{Co95} and Kawakita's classification of
extremal contractions to a singular point of the type under
consideration, see \cite{Kawakita2003}, --- there are precisely
two extremal contractions, which give, respectively, the two
varieties $Y$, mentioned above.

The set of polynomials $f\in {\cal P}_{4,5}$, such that the
quartic $\{f=0\}$ has a singularity of the type \cite{CoMe}, is of
codimension 4 in ${\cal P}_{4,5}$. A ``simpler'' or more typical
singularity of rank 3 (considered in the present paper) gives the
codimension 2 in ${\cal P}_{4,5}$ --- such quartics are much more
typical. However, their birational rigidity is stated in
\cite{CoMe} as a conjecture only. The classification of extremal
contractions to quadratic singularities of rank 3, given in
\cite{Kawakita2002}, does not lead to a proof of the birational
rigidity via an application of the Sarkisov program, and the
argument using the inversion of adjunction, mentioned above, does
not lead to a success, either: on the quadric cone $Q\subset
{\mathbb P}^3$ one can find infinitely many curves (that is,
effective 1-cycles) $C$, cut out by a surface of degree $m$, such
that the pair $(Q,\frac{1}{m} C)$ is not log canonical. For that
reason, the problem of description of the birational type of a
quartic with a double point of rank 3 has been open up to this
day.

However, the technique of the method of maximal singularities,
used in \cite{Pukh89c}, works successfully for such quartics,
either. In the present paper we replace the test class technique
by the technique of counting multiplicities (\cite{Pukh00c} or
\cite[Chapter 2]{Pukh13a}) and add the technical tool of modifying
the graph of a maximal singularity, which essentially simplifies
the arguments (in 1986, when \cite{Pukh89c} was written, this
technique has not yet been developed). Taking into account that
the birational rigidity is known for quartics with ordinary double
points, we obtain from Theorems 0.1 and 0.2, that the birational
type of three-dimensional quartics is now described in codimension
$\leqslant 2$.

To conclude, let us discuss briefly the following question: for
which codimensions it is realistic to expect birational rigidity
type results (where \cite{CoMe} belongs, too)? Quartics with a
triple point are rational and the set of their equations is of
codimension 11 in ${\cal P}_{4,5}$; quartics with a double line
have a structure of a conic bundle and the set of their equations
is of codimension 11, too; quartics containing a two-dimensional
plane have a structure of a fibration into cubic surfaces (for
what is known about their birational geometry, see
\cite{Grinenko2011}) and the set of their equations is of
codimension 9. It seems that one can expect the birational
rigidity type results in codimension $\leqslant 8$. See also the
paper \cite{AbbanKaloghiros2016}, belonging to this class of
problems, where a general conjecture on the birational rigidity of
three-dimensional quartics with isolated singularities is stated
that covers the claim of Theorem 0.1.\vspace{0.3cm}


{\bf 0.5. Acknowledgements.} The author is grateful to the members
of Divisions of Algebraic Geometry and Algebra at Steklov
Institute of Mathematics for the interest to his work, and also to
the colleagues-algebraic geometers at the University of Liverpool
for general support.


\section{Twenty five involutions}

In this section we reduce the proof of Theorems 0.1 and 0.2 to the
claim of Theorem 0.3, and the proof of the latter theorem to a
certain fact on the multiplicities of a mobile linear system on
$V$ (Theorem 1.1). The decisive tool in the arguments is the 25
birational involutions.\vspace{0.3cm}

{\bf 1.1. Factoriality.} First of all, let us show that the
assumptions of Subsection 0.1 imply that $V$ is a factorial
variety. Let $\sigma_{\mathbb P}\colon{\mathbb P}^+\to{\mathbb
P}^4$ be the blow up of the point $o$ and $E_{\mathbb P}$ its
exceptional divisor, which identifies naturally with the
projective space ${\mathbb P}^3_{(z_1:z_2:z_3:z_4)}$ in the
notations of Subsection 0.1. Let $V^+\subset{\mathbb P}^+$ be the
strict transform of $V$, so that
$$
\sigma=\sigma_{\mathbb P}|_{V^+}\colon V^+\to V
$$
is the blow up of the point $o$ on $V$ and $Q=V^+\cap E_{\mathbb
P}$ the exceptional quadric. It is easy to check that the
condition $q_3(0,0,0,1)\neq 0$ implies, that the variety $V^+$ is
non-singular. Let $H_{\mathbb P}$ be the class of a hyperplane in
${\mathbb P}^4$, then
$$
\mathop{\rm Pic}{\mathbb P}^+={\mathbb Z}H_{\mathbb
P}\oplus{\mathbb Z}E_{\mathbb P},
$$
where for simplicity we write $H_{\mathbb P}$ instead of
$\sigma^*_{\mathbb P}H_{\mathbb P}$. Now
$$
V^+\sim 4H_{\mathbb P}-2E_{\mathbb P}
$$
is a very ample non-singular divisor on ${\mathbb P}^+$, so that
by the Lefschetz theorem
$$
\mathop{\rm Pic}V^+={\mathbb Z}H\oplus{\mathbb Z}Q,
$$
which implies the factoriality of the quartic $V$.\vspace{0.3cm}


{\bf 1.2. Proof of Theorems 0.1 and 0.2.} Assuming the claim of
Theorem 0.3, let us prove our main results. Let $\pi'\colon V'\to
S'$ be a Mori fibre space, where $\mathop{\rm dim}S'\geqslant 1$.
For our argument it is sufficient to assume that the fibre of the
projection $\pi'$ over a point of general position on the base
$S'$ is either a rational curve (if $\mathop{\rm dim}S'=2$), or a
rational surface (if $\mathop{\rm dim}S'=1$). Assume that there is
a birational isomorphism $\chi\colon V\dashrightarrow V'$.

Fix an arbitrary very ample linear system on $S'$ and let
$\Sigma'$ be its pull back on $V'$, and
$\Sigma=(\chi^{-1})_*\Sigma'$ the strict transform of $\Sigma'$ on
$V$. It is well known that the pair $(V,\frac{1}{n}D)$, where
$D\in\Sigma$ is a general divisor and $n=n[\Sigma]$, is not
canonical. Applying Theorem 0.3, we see that, ``twisting'' $\chi$
by means of $\tau_a$, that is to say, replacing $\chi$ by
$\chi\circ\tau_a$, we obtain a new birational isomorphism
$V\dashrightarrow V'$, and moreover,
$n[(\tau_a)_*\Sigma]<n(\Sigma)$. This process does not terminate,
because the property of the pair $(V,\frac{1}{n}D)$ to be
non-canonical follows from the inequality $\mathop{\rm
dim}S'\geqslant 1$ and for that reason is preserved. We obtain an
infinite strictly decreasing sequence of positive integers. This
is impossible. Q.E.D. for the first claim of Theorem 0.1.

Let $V'$ be a Fano three-fold with ${\mathbb Q}$-factorial
singularities and the Picard number 1 and $\chi\colon
V\dashrightarrow V'$ a birational map. Let us fix a complete very
ample linear system $\Sigma'$ on $V'$ and set
$\Sigma=(\chi^{-1})_*\Sigma'$. Again applying the ``untwisting''
procedure, we may assume that already the pair $(V,\frac{1}{n}D)$,
where $D\in\Sigma$ is a general divisor and $n=n[\Sigma]$, is
canonical. In that case it is well known (see, for instance,
\cite[Chapter 2, Proposition 1.6]{Pukh13a}), that $\chi$ is a
biregular isomorphism. Q.E.D. for Theorem 0.1.

Let us show Theorem 0.2. Take $\chi\in\mathop{\rm Bir}V$, where
$\chi\not\in\mathop{\rm Aut}V$. For $\Sigma'$ we tale the complete
linear system of hyperplane sections of the quartic
$V\subset{\mathbb P}^4$. Obviously, $n[\Sigma]\geqslant 2$
(otherwise, $\chi$ is a biregular automorphism). The pair
$(V,\frac{1}{n}D)$ is not canonical (otherwise, as above, we would
have had $\chi\in\mathop{\rm Aut}V)$, so that we apply Theorem
0.3. Therefore, for some integer-valued sequence
$$
a(j)\in\{0,1,\dots,24\},
$$
$j=1,\dots,k$, we have
$$
\chi\circ\tau_{a(1)}\circ\dots\circ\tau_{a(k)}\in\mathop{\rm
Aut}V.
$$
Therefore, the subgroups $\mathop{\rm B}(V)$ and $\mathop{\rm
Aut}V$ generate $\mathop{\rm Bir}V$.

From the explicit construction of the involutions $\tau_i$, given
in Subsection 1.3, it is obvious that for every biregular
automorphism $\chi$
$$
\chi^{-1}\tau_0\chi=\tau_0,\quad\chi^{-1}\tau_i\chi=\tau_j
$$
for $i\geqslant 1$, where $\chi(L_j)=L_i$ (obviously, a projective
automorphism preserves the point $o$ and permutes the lines
$L_i$), so that $\mathop{\rm B}(V)$ is a normal subgroup.

Finally, let us show that all relations between the generators
$\tau_i$ follow from the the relations $\tau_i^2=\mathop{\rm
id}_V$ and that
$$
\mathop{\rm B}(V)\cap\mathop{\rm Aut}V=\{\mathop{\rm
id}\nolimits_V\}.
$$
Assume that there is a sequence $a(j)\in\{0,1,\dots,24\}$,
$j=1,\dots,k$, where $k\geqslant 2$, neither two neighbor integers
are equal, $a(j)\neq a(j+1)$ for $j=1,\dots,k-1$ and
$$
\tau_{a(1)}\circ\cdots\circ\tau_{a(k)}\in\mathop{\rm Aut}V.
$$
We may assume that this sequence is the shortest of all possible
ones, so that for $j\leqslant k-1$
$$
\tau_{a(1)}\circ\cdots\circ\tau_{a(j)}\not\in\mathop{\rm Aut}V
$$
Again, let $\Sigma'$ be the complete linear system of hyperplane
sections of the quartic $V$. Set for $j=1,\dots,k$
$$
\Sigma_j=(\tau_{a(j)}\circ\cdots\circ\tau_{a(1)})_*\Sigma',
$$
$\Sigma_j\subset|n_jH|$. For convenience set
$\Sigma_0=\Sigma'=|H|$ and $n_0=1$. Obviously, $n_k=1$, however
$n_j\geqslant 2$ for $j=1,\dots,k-1$. Clearly, $n_j=n[\Sigma_j]$
for $j=0,\dots,k$. Since
$$
\Sigma_j=(\tau_{a(j)})_*\Sigma_{j-1},
$$
Theorem 0.3 implies that $n_{j-1}\neq n_j$ for $j=1,\dots,k$. Let
$D_j\in\Sigma_j$ be a general divisor, then the pair
$$
\left(V,\frac{1}{n_j}D_j\right)
$$
is not canonical for $j=1,\dots,k-1$. Choose
$$
e\in\{1,\dots,k-1\}
$$
such that $n_e$ is a local maximum of the sequence $n_j$, that is,
$$
n_e>n_{e-1}\quad\mbox{and}\quad n_e>n_{e+1}.
$$
Since
$$
\Sigma_{e-1}=(\tau_{a(e)})_*\Sigma_e\quad \mbox{and}\quad
\Sigma_{e+1}=(\tau_{a(e+1)})_*\Sigma_e
$$
and $a(e)\neq a(e+1)$, we obtain a contradiction with Theorem 0.3.
This completes the proof of Theorems 0.1 and 0.2.\vspace{0.3cm}


{\bf 1.3. Birational involutions.} Recall the construction of
birational involutions $\tau_i$ \cite{Pukh89c}. Let
$L^+_1,\dots,L^+_{24}\subset V^+$ be the strict transforms of the
lines $L_1,\dots,L_{24}$ (see Subsection 1.1). Set
$$
\pi_0\colon V^+\to{\mathbb P}^3
$$
to be the (regularized) projection from the point $o$. Obviously,
$\pi_0$ is a map of degree 2 and its restriction onto
\begin{equation}\label{17.09.24.1}
V^+\backslash\left(\bigcup^{24}_{i=1}L^+_i\right)
\end{equation}
is a double cover of the complement
$$
{\mathbb P}^3\backslash\left(\bigcup^{24}_{i=1}\pi_0(L^+_i)\right)
$$
to the 24 points $\pi_0(L^+_i)\in{\mathbb P}^3$. The involution
$\tau_0$ is the Galois involution, corresponding to the map
$\pi_0$. Obviously, $\tau_0$ is a biregular involution of the open
set (\ref{17.09.24.1}).

{\bf Lemma 1.1.} {\it The action of the involution $\tau_0$ on
$\mathop{\rm Pic}V^+$ is given by the relations}
$$
\tau^*_0H=3H-4Q,\quad\tau^*_0Q=2H-3Q.
$$
(This is \cite[Lemma 1.2]{Pukh89c}.)

Now let us describe the involutions $\tau_i,1\leqslant i\leqslant
24$.

Let $V^+_i\to V^+$ be the blow up of the line $L^+_i$,
$\pi_i\colon V^+_i\to{\mathbb P}^2$ the (regularized) projection
from the line $L_i$. It is easy to see that for a general point
$t\in{\mathbb P}^2$ the fibre $C_t=\pi^{-1}_i(t)$ is an elliptic
curve, and moreover, $(C_t\cdot Q_i)=1$, where $Q_i\subset V^+_i$
is the strict transform of the quadric $Q$, so that the surface
$Q_i$ is a rational section of the fibration into elliptic curve.
We define $\tau_i$ as the fibre-wise reflection on $C_t$ with
respect to the point $C_t\cap Q_i$. Let $\Delta_i$ be the
exceptional divisor of the blow up of $L^+_i$, so that
$$
\mathop{\rm Pic}V^+_i={\mathbb Z}H\oplus{\mathbb Z}Q\oplus{\mathbb
Z}\Delta_i.
$$

{\bf Lemma 1.2.} {\it The birational involution $\tau_i$ extends
to a biregular involution of an open set, which is the complement
to a finite set of fibres of the projection $\pi_i$. Its action on
$\mathop{\rm Pic}V^+_i$ is given by the relations}
$$
\tau^*_iH=11H-6Q-12\Delta_i,\quad\tau^*_iQ=Q,\quad
\tau^*_i\Delta_i=10H-6Q-11\Delta_i
$$
(This is \cite[Lemma 1.3 and Remark 1.3]{Pukh89c}.)

Now we can start the {\bf proof of Theorem 0.3}. Fix a mobile
linear system $\Sigma$, set $n=n[\Sigma]$ and assume that the pair
(\ref{13.09.24.1}) is not canonical. Let $\Sigma^+$ be the strict
transform of $\Sigma$ on $V^+$,
$$
\Sigma^+\subset|nH-mQ|
$$
for some $m\in{\mathbb Z}_+$.

{\bf Theorem 1.1.} {\it Among the 25 integers
$$
m,\quad\mathop{\rm mult}\nolimits_{L_1} \Sigma,\quad \dots,\quad
\mathop{\rm mult}\nolimits_{L_{24}}\Sigma
$$
precisely one is strictly higher than $n$, and all the others are
strictly smaller than $n$.}

Assuming Theorem 1.1, let us show Theorem 0.3. By Lemma 1.1,
$$
n[(\tau_0)_*\Sigma]=3n-2m,
$$
and by Lemma 1.2 for $a\in\{1,\dots,24$\}
$$
n[(\tau_a)_*\Sigma]=11n-10\mathop{\rm mult}\nolimits_{L_a}\Sigma.
$$
Thus if $m>n$, then the inequality
$$
n[(\tau_0)_*\Sigma]<n
$$
holds, and moreover, $n[(\tau_i)_*\Sigma]>n$ for $i\geqslant 1$.
If $\mathop{\rm mult}_{L_a}\Sigma> n$, then the inequality
$$
n[(\tau_a)_*\Sigma]>n
$$
holds, and moreover, $[(\tau_i)_*\Sigma]>n$ for $i\neq a$. Proof
of Theorem 0.3 is complete.

The next two sections contain a proof of Theorem 1.1.


\section{The method of maximal singularities}

In this section we carry out some preparatory work for the proof
of Theorem 1.1. We show that the centre of a maximal singularity
$E$ on $V$ is either the point $o$, or one of the 24 lines $L_i$,
and in the latter case the claim of Theorem 1.1 holds. For the
case when the centre of $E$ is the point $o$, we state Theorem 2.1
that directly implies Theorem 1.1 and start to prove
it.\vspace{0.3cm}

{\bf 2.1. Non-singular points.} In the notations of Subsections
1.3 and 0.2 let $E$ be a maximal singularity of the linear system
$\Sigma$. Let us exclude the option (2). Assume that the centre of
$E$ on $V$ is a non-singular point $p\neq o$. Consider the {\it
self-intersection} $Z=(D_1\circ D_2)$ of the system $\Sigma$,
where $D_1,D_2\in\Sigma$ are general divisors and the symbol
$(D_1\circ D_2)$ means the effective 1-cycle of the
scheme-theoretic intersection of $D_1$ and $D_2$. Since $\Sigma$
has no fixed components, the self-intersection is well defined. It
is well known that the $4n^2$-inequality holds:
$$
\mathop{\rm mult}\nolimits_pZ>4n^2.
$$
At the same time, the effective cycle $Z$ as a 1-cycle on
${\mathbb P}^4$ is of degree $4n^2$. We get a contradiction,
excluding the option (2). This argument goes back to the classical
paper \cite{IM}, where the language of the {\it test class} is
used, see \cite[Chapter 2]{Pukh13a}. Various proofs of the
$4n^2$-inequality and its comparison with the test class technique
see in \cite[Chapter 2, \S 2 and \S 3]{Pukh13a}.\vspace{0.3cm}


{\bf 2.2. Curves.} Assume now that the centre of the maximal
singularity $E$ is an irreducible curve $B\subset V$. In that case
it is easy to check that the inequality
$$
\mathop{\rm mult}\nolimits_B\Sigma>n
$$
holds, that is, $B$ is a {\it maximal curve} of the system
$\Sigma$, see \cite[Chapter 2]{Pukh13a}. We already saw in \S 1,
that there are mobile linear systems with one of the 24 lines
$L_i$ as a maximal curve. Assume that $B$ is not a line passing
through the point $o$. Considering again the self-intersection $Z$
and taking into account that $B$ comes into $Z$ with the
multiplicity $> n^2$, we see that $\mathop{\rm deg}B\leqslant 3$.
If $o\not\in B$, then we get a contradiction, arguing in the word
for word the same way as in \cite{IM}. Assume therefore that $o\in
B$.

If the curve $B$ is not contained in a plane, that is, the linear
span $\langle B\rangle$ has dimension $\geqslant 3$, then $B$ is a
rational normal curve and $\mathop{\rm dim}\langle B\rangle=3$.
Set $\mu=\mathop{\rm mult}_B\Sigma$ and let $B^+\subset V^+$ be
the strict transform of $B$ on $V^+$. The curve $B$ is cut out by
quadrics, so that consider the surface $S$, which is the
intersection of $V$ with a general quadric containing $B$.
Obviously, $S$ is non-singular outside the point $o$. Let $S^+$ be
its strict transform on $V^+$.

{\bf Lemma 2.1.}  {\it The surface $S^+$ is non-singular.}

{\bf Proof.} This is almost obvious. If the point $p=B^+\cap Q$ is
not the vertex of the quadric cone $Q$, then $S^+\cap Q$ is a
non-singular point, so that $S^+$ is non-singular. If the point
$p$ is the vertex of the cone $Q$, then $p$ \, is the only
candidate for a singularity of the surface $S^+$, because $S^+\cap
Q$ is a pair of distinct lines, intersecting at the point $p$.
However, $V^+$ is non-singular, so that the tangent hyperplane
$$
T_pV^+\subset T_p{\mathbb P}^+
$$
coincides with the tangent hyperplane $T_pE_{\mathbb P}$, whereas
the quadric hypersurface that cuts out $S$ on $V$ is non-singular
at the point $o$, and for that reason the tangent hyperplane to
its strict transform on ${\mathbb P}^+$ at the point $p$ can not
be equal to $T_pE_{\mathbb P}$. Therefore, $S^+$ is non-singular
in any case. Q.E.D. for the lemma.

The restriction $\Sigma^+|_{S^+}$ is a linear system of curves on
the non-singular surface $S^+$ with the non-singular rational
curve $B^+$ as its only fixed component. Computing the
self-intersection of the mobile part of that system, we get the
inequality
$$
8n^2-6n\mu-m^2-4\mu^2-(\mu-m)^2\geqslant 0,
$$
see \cite[Subsection 4.3]{Pukh89c}. This is impossible for
$\mu>n$. Therefore, the curve $B$ is contained in a
plane.\vspace{0.3cm}


{\bf 2.3. The plane section lemma.} Let $P\ni o$ be a plane in
${\mathbb P}^4$. The plane curve $V\cap P$ is a union of
irreducible components $C_1,\dots,C_k$ of degrees $d_1,\dots,d_k$,
where $k\leqslant 4$.

{\bf Proposition 2.1 (the plane section lemma).} (i) {\it If
$d_i\geqslant 2$, then} $\mathop{\rm mult}_{C_i}\Sigma\leqslant
n$.

(ii) {\it If $d_i=d_j=1,i\neq j$, then}
$$
\mathop{\rm mult}\nolimits_{C_i}\Sigma+\mathop{\rm
mult}\nolimits_{C_j}\Sigma\leqslant 2n.
$$

{\bf Proof} is word for word the same as the arguments in \cite[\S
6]{Pukh89c}, proving the plane section lemma in \cite[Subsection
4.2]{Pukh89c}; the fact that the rank of the quadratic form
$q_2(z_*)$ (see Subsection 0.1) is equal to 4, is never used in
\cite[\S 6]{Pukh89c}; only the irreducibility of the exceptional
quadric $\sigma^{-1}(o)$ of the blow up of the point $o$ is used,
and in the case under consideration this is true. Q.E.D. for the
proposition.

The plane section lemma shows that the option (3) (Subsection 0.2)
does not take place, and if $\mathop{\rm mult}_{L_i}\Sigma>n$
(that is, the option (1) takes place), then for $j\neq i$ we have
$\mathop{\rm mult}_{L_j}\Sigma<n$ (every pair of lines $L_i,L_j$
is contained in some plane).\vspace{0.3cm}


{\bf 2.4. The uniqueness in Theorem 1.1.} Let $L=L_i$, $1\leqslant
i\leqslant 24$, be one of the lines. Set $\mu=\mathop{\rm
mult}_{L}\Sigma$. The integer $m=\mathop{\rm
ord}\nolimits_Q\Sigma$ was defined in Subsection 1.3.

{\bf Lemma 2.2.} {\it If one of the numbers $m,\mu$ is strictly
higher than $n$, then the other one is strictly smaller than $n$.}

{\bf Proof.} We may assume that $\mu\geqslant n$. Let $S$ be the
section of the quartic $V$ by a general hyperplane, containing the
line $L$, and $S^+$ its strict transform on $V^+$. Obviously,
$S^+$ is a non-singular surface and the restriction
$\Sigma^+|_{S^+}$ is a linear system of curves with the only fixed
component $L^+$ of multiplicity $\mu$. It is easy to compute (see
\cite[Subsection 4.4.]{Pukh89c}), that the self-intersection of
the mobile part of the system $\Sigma^+|_{S^+}$ (after removing
$L^+$) is
$$
4n^2-m^2-\mu^2-2n\mu-(m-\mu)^2\geqslant 0,
$$
whence we get $m\leqslant n$, which proves the lemma.

Thus if one of the 25 integers in the statement of Theorem 1.1 is
strictly higher than $n$, then all the rest are strictly smaller
than $n$. It remains to show (and this is the hardest part of the
work), that at least one of the 25 integers is strictly higher
than $n$. If the option (1) takes place (see Subsection 0.2), then
the inequality $\mathop{\rm mult}_{L_i}> n$ holds, where $L_i$ is
the centre of the maximal singularity $E$ on $V$, and there is
nothing to prove.

Therefore, starting from this point, we may assume that the option
(0) takes place, that is, the centre of (any, if there are more
than one) maximal singularity $E$ is the point $o$.\vspace{0.3cm}


{\bf 2.5. The plan of further work.} We will prove the following
claim (which immediately implies Theorem 1.1).

{\bf Theorem 2.1.} {\it Assume that the centre of the maximal
singularity $E$ on $V$ is the point $o$, and for $i=1,\dots,24$
$$
\mathop{\rm mult}\nolimits_{L_i}\Sigma\leqslant n.
$$
Then the inequality $m>n$ holds (that is, already $Q$ is a maximal
singularity of the system $\Sigma$).}

{\bf Proof}  starts here but its main part makes the contents of
\S 3. Assume that $m\leqslant n$. Then the centre of $E$ on $V^+$
is an irreducible closed subset of the cone $Q$.

{\bf Lemma 2.3.} {\it The centre of $E$ on $V^+$ is not a curve of
degree $\geqslant 2$ in $E_{\mathbb P}$}.

{\bf Proof.} Assume the converse: the centre of $E$ on $V^+$ is an
irreducible curve $C\subset Q$ of degree $d_C\geqslant 2$. The
restriction $\Sigma_Q$ of the mobile system $\Sigma^+$ onto $Q$ is
cut out on $Q$ by surfaces of degree $m\leqslant n$ in $E_{\mathbb
P}$, so that the curves of that non-empty system are of degree
$2m$. However, the inequality
$$
\mathop{\rm mult}\nolimits_C\Sigma^+> n
$$
holds, so that $\Sigma_Q$ must contain the curve $C$ of degree
$d_C\geqslant 2$ with multiplicity $>n$. This is impossible.
Q.E.D. for the lemma.

{\bf Lemma 2.4.} {\it The centre of $E$ on $V^+$ is not the vertex
$v\in Q$ of the cone $Q$}.

{\bf Proof.} Assume the converse: the centre of $E$ on $V^+$ is
the point $v$. Let $D\in\Sigma$ be a general divisor,
$D^+\in\Sigma^+$ its strict transform on $V^+$, $D^+\sim nH-mQ$.
Since $m\leqslant n$, the inequality
$$
\mathop{\rm ord}\nolimits_ED^+>n\cdot a(E,V^+)
$$
holds, which implies that
$$
\mathop{\rm mult}\nolimits_vD^+> n.
$$
Recall that the variety $V^+$ is non-singular. Let $R\subset Q$ be
an arbitrary line (a generator of the cone $Q$). By what was said
above we get the inequality
$$
(R\cdot D^+)_v> n.
$$
At the same time,
$$
(R\cdot D^+)=m(R\cdot (-Q))=m(R\cdot (-E_{\mathbb P}))_{{\mathbb
P}^+}=m\leqslant n.
$$
This is only possible if $R\subset D^+$. However, the generators
sweep out the cone $Q$ and $Q\not\subset D^+$. This contradiction
proves Lemma 2.4.

{\bf Lemma 2.5.} {\it The centre of $E$ on $V^+$ is not a line.}

{\bf Proof.} Assume the converse: the centre of $E$ is a generator
$G\subset Q$; then $\mathop{\rm mult}_G\Sigma^+> n$ and, in
particular, $\mathop{\rm mult}_v\Sigma^+> n$, where $v\in Q$ is
the vertex of the cone. Now we argue as in the proof of Lemma 2.4.
Q.E.D. for the lemma.

We conclude that the centre of the maximal singularity $E$ on
$V^+$ is a point $p\in Q$, which is not the vertex of the cone
$Q$.

Our further work is organized as follows: we show

\noindent --- that the point $p$ lies on the strict transform
$L^+_a$ of one of the 24 lines $L_1,\dots,L_{24}$ on $V$,
$p=L^+_a\cap Q$,

\noindent --- that the input of the line $B=L_a$ into the
self-intersection $Z=(D_1\circ D_2)$ of the system $\Sigma$ is
quite large,

\noindent --- and that the latter is incompatible with the
assumptions of Theorem 2.1.

All this will be done in \S 3.


\section{The technique of counting multiplicities}

This section contains the proof of Theorem 2.1. The arguments
follow the same scheme as in \cite[\S 5]{Pukh89c}, however with
the test class technique replaced by the technique of counting
multiplicities. Besides, we use the technical tool of modifying
the graph of the maximal singularity, which makes it possible to
work without the ``graph lemma'' \cite[\S 7]{Pukh89c}, simplifying
the proof.\vspace{0.3cm}


{\bf 3.1. Resolution of the maximal singularity.} following the
standard procedure of the method of maximal singularities
\cite[Chapter 2]{Pukh13a}, consider the sequence of blow ups
$$
\varphi_{i,i-1}\colon V_i\to V_{i-1},
$$
$i=1,\dots,K$, which is defined in the following way. The variety
$V_0$ is $V$, the blow up $\varphi_{1,0}$ is the blow up of the
point $o$, so that $V_1\cong V^+$. In general, $\varphi_{i,i-1}$
blows up the centre $B_{i-1}$ of the singularity $E$ on $V_{i-1}$,
the exceptional divisor $\varphi^{-1}_{i,i-1}(B_{i-1})$ is denoted
by the symbol $E_i$ and the centre of $E$ on $V_K$ is the
exceptional divisor $E_K$. Setting for $i>j$
$$
\varphi_{i,j}=\varphi_{j+1,j}\circ\dots\circ\varphi_{i,i-1}\colon
V_i\to V_j,
$$
we can say that $\varphi^{-1}\circ\varphi_{K,1}\colon
V_K\dashrightarrow\widetilde{V}$ (in the notations of Subsection
0.2) maps $E_K$ onto $E$ (this map is well defined at the general
point of $E_K$). Identifying $V_1$ with $V^+$, we identify $E_1$
with $Q$ in the previous notations. As we proved at the end of \S
2, the centre $B_1$ is a point $p\in Q$, which is not the vertex
of the cone $Q$.

By construction, $\varphi_{i,i-1}(B_i)=B_{i-1}$, so that the
dimensions $\mathop{\rm dim}B_i$ form a non-decreasing sequence
and we can break the set $\{1,\dots,K\}$ into the {\it lower}
$\{1,\dots,L\}$ and the {\it upper} $\{L+1,\dots,K\}$ parts: for
the ``lower'' indices $i\leqslant L$ the centre of the blow up
$B_{i-1}$ is a point, for the ``upper'' ones it is a curve. Let
$$
\delta_1=1,\quad \delta_2=\dots=\delta_L=2,\quad
\delta_{L+1}=\dots=\delta_K=1
$$
be the elementary discrepancies. On the set $\{1,\dots,K\}$ there
is a natural structure of an oriented graph $\Gamma$: the vertices
$i$ and $j$ are connected by an oriented edge (arrow), if and only
if $i>j$ and $B_{i-1}\subset E^{i-1}_j$, where the upper index
$i-1$ means the strict transform on $V_{i-1}$; in this case we
write $i\to j$. For $i>j$ the number of paths in $\Gamma$ from $i$
to $j$ is denoted by the symbol $p_{ij}$ and for convenience we
set $p_{ii}=1$.

Set $\nu_1=m$ (see Subsection 1.3) and for $i\geqslant 2$
$$
\nu_i=\mathop{\rm mult}\nolimits_{B_{i-1}}\Sigma^{i-1}.
$$
Since for a general divisor $D\in\Sigma$ the 1-cycle $(D^+\circ
Q)=(D^1\circ E_1)$ is an effective cycle of degree $2m=2\nu_1$, we
obtain the inequality $2\nu_1\geqslant\nu_2$. Besides, for
$i\geqslant 2$
$$
\nu_i\geqslant\nu_{i+1}.
$$
(See \cite[Chapter 2, \S\S 1-2]{Pukh13a} for the definition and
elementary properties of the graph $\Gamma$.) Thus for $i\geqslant
2$ we have $\nu_i\leqslant 2n$.\vspace{0.3cm}


{\bf 3.2. The Noether-Fano inequality.} The inequality
(\ref{21.09.24.1}) of Subsection 0.2 now can be re-written in the
following traditional form:
\begin{equation}\label{21.09.24.2}
\sum^K_{i=1}p_{Ki}\nu_i>n\sum^K_{i=1}p_{Ki}\delta_i.
\end{equation}
Every path from the vertex $K$ starts with an arrow $K\to a$,
where $a\leqslant K-1$. If $\nu_K\leqslant n$, re-write then the
inequality (\ref{21.09.24.2}) in the form
$$
(\nu_K-n)+\sum_{K\to
a}\left(\sum^a_{i=1}p_{ai}(\nu_i-n\delta_i)\right)>0,
$$
which implies that already some exceptional divisor $E_a$,
$a\leqslant K-1$, is a maximal singularity of the system $\Sigma$.
If we still have $\nu_a\leqslant n$, then this procedure can be
applied to $E_a$. Repeating, we obtain in the end a maximal
singularity $E_j$ such that $\nu_j>n$. In order to simplify the
notations, we assume that the inequality $\nu_K> n$ holds already
for $E$.

It is easy to see that the ``upper'' part of the graph $\Gamma$ is
a chain:
$$
(L+1)\leftarrow\dots\leftarrow K,
$$
that is, there are no other arrows between these vertices apart
from the subsequent ones $(i+1)\to i$, as in the case of an arrow
$(i+2)\to i$ for $i\geqslant L+1$ we would have
$$
\nu_i\geqslant\nu_{i+1}+\nu_{i+2}>2n,
$$
contrary to what we know; one should also take into account that
from $i\to j$ for $i\geqslant j+2$ it follows that $a\to j$ for
$j<a<i-1$ (this is true for the whole graph $\Gamma$, see
\cite[Chapter 2, \S\S 1-2]{Pukh13a}). It is also easy to check
that $(L+2)\nrightarrow L$. However, it is quite possible to have
arrows going from the ``upper'' vertices $L+1,\dots,K$ to the
``lower'' ones $1,\dots, L$, which are different from $(L+1)\to
L$.

If this is the case, we apply the known procedure of modifying the
graph $\Gamma$, removing (``erasing'') the arrows, joining the
vertices $L+1,\dots,K$ with the vertices $1,\dots,L-1$. We obtain
a new graph $\Gamma^+$ with the same set of vertices, but fewer
arrows. Note that in $\Gamma^+$ between the lower vertices
$1,\dots,L$ and between the upper vertices $L+1,\dots,K$ the
arrows are the same as in $\Gamma$. Let $p^+_{ij}$ be the number
of paths in $\Gamma^+$ from $i$ to $j$, if $i>j$ and $p^+_{ii}=1$.
Once again, note that $p^+_{ij}=p_{ij}$, if $j\leqslant i\leqslant
L$ or $L+1\leqslant j\leqslant i$.

{\bf Lemma 3.1.} {\it The Noether-Fano type inequality holds:}
\begin{equation}\label{23.09.24.1}
\sum^K_{i=1}p^+_{Ki}\nu_i>n\sum^K_{i=1}p^+_{Ki}\delta_i.
\end{equation}

{\bf Proof.} Re-write the inequality (\ref{21.09.24.2}) in the
following form:
$$
\sum^L_{i=1}p_{Ki}(\nu_i-n\delta_i)+\sum^K_{i=L+1}p_{Ki}(\nu_i-n\delta_i)>0.
$$
In the first sum all expressions $\nu_i-n\delta_i$ are
non-positive, in the second sum they are positive. Taking into
account that $p_{Ki}=p^+_{Ki}$ for $i\geqslant L+1$ and
$p^+_{Ki}\leqslant p_{Ki}$ for $i\leqslant L$, replacing $p_{Ki}$
by $p^+_{Ki}$, we can only make the Noether-Fano inequality
stronger. Q.E.D. for the lemma.

Set $r_i=p^+_{Ki}$. Obviously,
$r_L=r_{L+1}=\dots=r_K=1$.\vspace{0.3cm}


{\bf 3.3. The technique of counting multiplicities.} Let
$Z=(D_1\circ D_2)$ be the self-intersection of the linear system
$\Sigma$, where $D_1,D_2\in\Sigma$ are general divisors. The
symbol $Z^j$ means the strict transform of $Z$ on $V_j$. Set
$$
m_i=\mathop{\rm mult}\nolimits_{B_{i-1}}Z^{i-1}
$$
for $i=1,\dots,L$ (the centres $B_{i-1}$ of blow ups are points);
in particular, $m_1=\mathop{\rm mult}_oZ$ and $m_2=\mathop{\rm
mult}_pZ^+$. Obviously, the multiplicities $m_i$ do not increase:
$m_1\geqslant m_2\geqslant\dots\geqslant m_L$. The technique of
counting multiplicities give us an estimate from below for some
positive linear combination of the multiplicities $m_1,\dots,m_L$,
see \cite[Chapter 2, \S 2]{Pukh13a}.

Recall that a function
$$
a\colon\{1,\dots,L\}\to{\mathbb R}_+
$$
{\it is compatible with the structure of the graph} $\Gamma$, if
for any $i$, $1\leqslant i\leqslant L-1$,
$$
a(i)\geqslant\sum_{j\to i}a(j)
$$
(since $a(\cdot)$ is defined for $i\leqslant L$, the summation is
taken over all $j\leqslant L$, such that $j\to i$). Now the main
result of the technique of counting multiplicities is stated in
the following way.

{\bf Proposition 3.1.} {\it The following inequality holds:}
$$
\sum^L_{i=1}a(i)m_i\geqslant2a(1)\nu^2_1+
\sum^L_{i=2}a(i)\nu^2_i+a(L)\sum^K_{i=L+1}\nu^2_i.
$$

{\bf Proof} repeats word for word the arguments in the proof
\cite[Chapter 2, Proposition 2.4]{Pukh13a} with the only change:
since the first exceptional divisor $E_1=Q$ is an irreducible
quadric, the first equality on the top of \cite[p. 53]{Pukh13a}
should be replaced by
$$
2\nu^2_1+d_1=m_{0,1},
$$
where $m_{0,1}=m_1=\mathop{\rm mult}_oZ$. Q.E.D. for the
proposition.

It is easy to see that the function $a(i)=r_i$ is compatible with
the structure of the graph $\Gamma$: since in $\Gamma^+$ there are
no arrows, apart from $(L+1)\to L$, that go from the upper part to
the lower part, we have the following equality for $i\leqslant
L-1$:
$$
r_i=\sum_{L\geqslant j\to i}r_j.
$$
At the same time, for the indices of the lower part the subgraphs
of $\Gamma$ and $\Gamma^+$ are the same. Since $r_i=1$ for
$i=L,L+1,\dots,K$, we have just proved the following claim.

{\bf Proposition 3.2.} {\it The following inequality holds:}
$$
\sum^L_{i=1}r_im_i\geqslant 2r_1\nu^2_1+\sum^K_{i=2}r_i\nu^2_i.
$$

Minimizing the quadratic form in the right hand side in
$\nu_1,\dots,\nu_K$ on the hyperplane, the equation of which we
obtain, replacing in the inequality (\ref{23.09.24.1}) the sign
$>$ by the equality sign, we see that the minimum is attained for
$$
2\nu_1=\nu_2=\dots=\nu_K.
$$
Set $\Sigma_0=\sum\limits^L_{i=2}r_i$ and
$\Sigma_1=\sum\limits^K_{i=L+1}r_i=K-L$.

Now the obvious computations give the proof of the following fact.

{\bf Proposition 3.3.} {\it The following inequality holds:}
\begin{equation}\label{23.09.24.2}
\sum^L_{i=1}r_im_i\geqslant
2\frac{(r_1+2\Sigma_0+\Sigma_1)^2}{r_1+2\Sigma_0+2\Sigma_1}n^2.
\end{equation}

{\bf Corollary 3.1.} {\it The following inequality holds:}
$$
m_1+m_2>6n^2.
$$

{\bf Proof.} The inequality (\ref{23.09.24.2}) only gets stronger
if in the left hand part we replace $m_i$ for $i\geqslant 2$ by
$m_2$. Therefore, we get the estimate
$$
r_1m_1+\Sigma_0m_2>
2\frac{(r_1+2\Sigma_0+\Sigma_1)^2}{r_1+2\Sigma_0+2\Sigma_1}n^2.
$$
By the construction of the graph $\Gamma^+$, we get the estimate
$r_1\leqslant\Sigma_0$; besides, $m_1\geqslant m_2$. Taking into
account these facts, we see that the minimum of the sum
$(m_1+m_2)$ with the value $r_1m_1+\Sigma_0m_2$ fixed is attained
for $m_1=m_2$, so that the claim of the corollary follows from the
estimate
$$
2\frac{(r_1+2\Sigma_0+\Sigma_1)^2}{(r_1+\Sigma_0)(r_1+
2\Sigma_0+2\Sigma_1)}\geqslant 3,
$$
which is equivalent to the inequality 
$$
0\geqslant
r^2_1+r_1(\Sigma_0+2\Sigma_1)-(2\Sigma^2_0+2\Sigma_0\Sigma_1+2\Sigma^2_1).
$$
Here the right hand part is an increasing function of $r_1$, so
that it is sufficient to set in the right hand side $r_1=\Sigma_0$
and get $-2\Sigma^2_1$. Q.E.D. for the corollary.\vspace{0.3cm}


{\bf 3.4. The existence of a line.} Consider the cycle $Z$ as a
1-cycle on ${\mathbb P}^4$. Let $B\subset{\mathbb P}^4$ be the
line, passing through the point $o$ in the direction of the
infinitely near point $p$, that is, $B^+\ni p$ (so that $p=B^+\cap
E_{\mathbb P})$.

{\bf Proposition 3.4.} {\it The line $B$ is a component of the
cycle $Z$ of multiplicity $\beta>2n^2$. In particular, $B\subset
V$ is one of the 24 lines} $L_1,\dots,L_{24}$.

{\bf Proof.} Consider the linear system $|H_{\mathbb P}-B|$ of
hyperplanes in ${\mathbb P}^4$, containing the line $B$. Let
$R\in|H_{\mathbb P}-B|$ be a general divisor in this linear
system, $R^+$ its strict transform on ${\mathbb P}^+$. Obviously,
$o\in R$ and $p\in R^+$, so that for any irreducible curve $C\neq
B$ in ${\mathbb P}^4$ we have
$$
\mathop{\rm deg}C=(C\cdot R)\geqslant(C\cdot
R)_o\geqslant\mathop{\rm mult}\nolimits_oC+\mathop{\rm
mult}\nolimits_pC^+.
$$
If the line $B$ is not a component of the cycle $Z$, then the
inequality
$$
\mathop{\rm deg}Z=4n^2\geqslant m_1+m_2>6n^2
$$
holds, which is impossible. Therefore, $B$ is a component of $Z$
and, therefore, $B\subset V$ is one of the 24 lines on $V$,
passing through the point $o$.

Write $Z=Z_*+\beta B$, where the effective 1-cycle $Z_*$ does not
contain $B$ as a component. Set
$$
m^*_1=\mathop{\rm mult}\nolimits_oZ_*,\quad m^*_2=\mathop{\rm
mult}\nolimits_pZ^+_*.
$$
By what was said above,
$$
\mathop{\rm deg}Z_*=4n^2-\beta\geqslant m^*_1+m^*_2>6n^2-2\beta,
$$
whence we get $\beta>2n^2$. Q.E.D. for the
proposition.\vspace{0.3cm}


{\bf 3.5. Blowing up the line $B$.} Let
$$
\varphi_B\colon V_B\to V^+
$$
be the blow up of the curve $B^+$ and $E_B\subset V_B$ the
exceptional divisor. By the symbol $|H-B|$ we denote the linear
system of hyperplane sections, containing the line $B$ (so that
$|H-B|=|H_{\mathbb P}-B|_V$ in the sense of Subsection 3.4) and
let $\Theta$ be its strict transform on $V_B$.

{\bf Proposition 3.5.} (i) $E_B\cong{\mathbb P}^1\times{\mathbb
P}^1$ {\it and $\varphi_B\colon E_B\to B^+$ is the projection onto
the direct factor.}

(ii) {\it The linear system $\Theta$ is free and its restriction
onto $E_B$ is a free linear system of curves of bi-degree (1,1).}

{\bf Proof} is elementary, it repeats word for word the proof of
\cite[Lemma 1.1]{Pukh89c}. The idea is that the blow up
$\varphi_B$ can be looked at as the restriction of the blow up
$$
\varphi_{{\mathbb P},B}\colon{\mathbb P}_B\to{\mathbb P}^+
$$
of the curve $B^+$ on ${\mathbb P}^+$ onto the subvariety $V^+$.
its exceptional divisor ${\mathbb E}_B\subset{\mathbb P}_B$ is
$B^+\times{\mathbb P}^2=B\times{\mathbb P}^2$. Considering the
normal sheaf ${\cal N}_{B^+\slash V^+}$, from the exact sequence
$$
0\to{\cal N}_{B^+\slash S^+}\to{\cal N}_{B^+\slash V^+}\to{\cal
O}(S^+)|_{B^+}\to 0,
$$
where $S\in|H-B|$ is a general surface and $S^+$ its strict
transform on $V^+$, we conclude that for ${\cal N}_{B^+}\slash
V^+$ there are two options: either ${\cal O}_{B^+}(-1)\oplus{\cal
O}_{B^+}(-1)$, or ${\cal O}_{B^+}(-2)\oplus{\cal O}_{B^+}$ (where
$B^+\cong{\mathbb P}^1$). In the first case we get the claim (i).
In the second case $E_B$ is a linear surface of type ${\mathbb
F}_2$, however, it is easy to see that the exceptional section of
that surface is a curve of the type
$$
B\times\mbox{\{point\}}
$$
in the sense of the identification of $B\times{\mathbb P}^2$ and
${\mathbb E}_B$. However, the presence of such curve on $V_B$
means that the scheme-theoretic intersection of surfaces
$\{q_2=0\}$, $\{q_3=0\}$, $\{q_4=0\}$ in ${\mathbb P}^3$ (see
Subsection 0.1) at the point, corresponding to the line $B$, has
multiplicity $\geqslant 2$, that is, there are $\leqslant 23$
distinct lines through the point $o$ on the quartic $V$, which
contradicts our assumption. The claim (ii) is obvious (see
\cite[Lemma 1.1]{Pukh89c}). Q.E.D. for the
proposition.\vspace{0.3cm}


{\bf 3.6. Restriction onto the surface $S$.} Again let us consider
a general surface $S\in|H-B|$. Its strict transform $S^+\subset
V^+$ is non-singular and the strict transform $S^*\subset V_B$
intersects the exceptional divisor $E_B$ by a non-singular curve
of bi-degree (1,1). Set
$$
\nu_B=\mathop{\rm mult}\nolimits_B\Sigma.
$$
Obviously, the linear system $\Sigma^+|_{S^+}$ is a linear system
of curves with a unique fixed component $B^+$ of multiplicity
$\nu_B$. By the symbol $\Sigma^B$ denote the strict transform of
the system $\Sigma$ on $V_B$, so that if $D_1,D_2\in\Sigma$ are
general divisors, then their strict transforms
$D^B_1,D^B_2\in\Sigma^B$ are general divisors of that system. Let
$Z_B$ be the part of the effective 1-cycle $(D^B_1\circ D^B_2)$
(the self-intersection of the system $\Sigma^B$), the support of
which is contained in the exceptional divisor $E_B\cong{\mathbb
P}^1\times B^+$. Let $f_B\in\mathop{\rm Pic}E_B$ be the class of a
fibre of the projection $\varphi_B|_{E_B}$ on $B^+$.

We can consider the effective 1-cycle $Z_B$ as an effective
divisor on $E_B$. Set
$$
d_B=(f_B\cdot Z_B).
$$
Thus $Z_B$ is a divisor of bi-degree $(d_B,*)$.

{\bf Lemma 3.2.} {\it The following equality holds:}
$$
\beta=\nu_B^2+d_B.
$$

{\bf Proof:} this is obvious.

The surfaces $S^+$ and $S^*$ are isomorphic: $S^*$ is the blow up
of the non-singular curve $B^+$ on the non-singular surface $S^+$.
The restriction of the linear system $\Sigma^B$ onto $S^*$ has no
fixed components. Since $\mathop{\rm mult}\nolimits_p
\Sigma^+=\nu_2$ and $p\in B^+$, the linear system $\Sigma^B$ has
the multiplicity $(\nu_2-\nu_B)$ along the fibre
$\varphi^{-1}_B(p)$, so that the linear system of curves
$\Sigma^B|_{S^*}$ has the base point $S^*\cap \varphi^{-1}_B(p)$
of multiplicity $(\nu_2-\nu_B)$. Finally, every irreducible curve
in the 1-cycle $Z_B$ that covers the base $B^+$ meets $S^*$ at
$\geqslant(f_B\cdot C)$ points of general position on $C$, which
are different from the point $S^*\cap\varphi^{-1}_B(p)$. It is
easy to compute that the self-intersection of the mobile system of
curves $\Sigma^B|_{S^*}$ is equal to
$$
4n^2-2m^2+2m\nu_B-2n\nu_B-2\nu^2_B=
4n^2-\nu^2_1-(\nu_1-\nu_B)^2-2n\nu_B-\nu^2_B,
$$
so that, taking into account what was said above, we obtain the
following claim.

{\bf Proposition 3.6.} {\it The following inequality holds:}
\begin{equation}\label{25.09.24.1}
4n^2-\nu^2_1-(\nu_1-\nu_B)^2-(\nu_2-\nu_B)^2-2n\nu_B\geqslant\beta.
\end{equation}\vspace{0.3cm}


{\bf 3.7. End of the proof of Theorem 2.1.} It is easy to check
that for $\nu_B\leqslant n$ the maximum of the right hand side of
the inequality (\ref{25.09.24.1}) on the ray
$\nu_1+\nu_2=3\theta$, $\nu_2\leqslant 2\nu_1$, is attained for
$\nu_1=\theta$ and $\nu_2=2\theta$ and equal to $\Lambda(\theta)$,
where
$$
\Lambda(t)=4n^2-2n\nu_B-t^2-(t-\nu_B)^2-(2t-\nu_B)^2.
$$

{\bf Lemma 3.3.} {\it For $t>\frac{n}{2}$ the function
$\Lambda(t)$ is decreasing.}

{\bf Proof:} $\Lambda'(t)=6(-2t+\nu_B)$ and by assumption
$\nu_B\leqslant n$. Q.E.D. for the lemma.

Furthermore, it is easy to check that the maximum of the value
$\Lambda(\theta)$ as a function of $\nu_B$ is attained for
$\nu_B=\frac12(3\theta-n)\in\left(\frac{n}{4},n\right)$ and equal
to
$$
\frac{9}{2}n^2-3n\theta-\frac{3}{2}\theta^2.
$$
Now it is easy to obtain from the inequality
$\Lambda(\theta)>2n^2$ that
$$
\theta<\frac{2\sqrt{6}-3}{3}n<\frac23n.
$$
On the other hand, from the Noether-Fano type inequality
(\ref{23.09.24.1}) we get:
$$
r_1\nu_1+(\Sigma_0+\Sigma_1)\nu_2>n(r_1+2\Sigma_0+\Sigma_1),
$$
and the maximum of the left hand side on the ray
$\nu_1+\nu_2=3\theta$, $\nu_2\leqslant 2\nu_1$, is attained again
for $\nu_1=\theta$ and $\nu_2=2\theta$, whence we get the estimate
$$
\Sigma_1>\frac{n-\theta}{2\theta-n}(r_1+2\Sigma_0).
$$
As we saw above, $\theta<\frac{2}{3}n$, so that finally
$$
\Sigma_1>r_1+2\Sigma_0.
$$

{\bf Proposition 3.7.} {\it The following inequality holds:}
$$
\frac{(r_1+2\Sigma_0+\Sigma_1)^2}{(r_1+\Sigma_0)(r_1+2\Sigma_0+
2\Sigma_1)}\geqslant 2.
$$

{\bf Proof.} This inequality is equivalent to the inequality
$$
0\geqslant r^2_1+2r_1\Sigma_0+2r_1\Sigma_1-\Sigma^2_1.
$$
The derivative of the right hand side by $\Sigma_1$ is
$2(r_1-\Sigma_1)<0$, so that it is sufficient to replace
$\Sigma_1$ by $r_1+2\Sigma_0$ and show the inequality
$$
0\geqslant 2r^2_1+2r_1\Sigma_0-4\Sigma^2_0,
$$
which is true since $r_1\leqslant\Sigma_0$. Q.E.D. for the
proposition.

Now from Proposition 3.3 we get
$$
r_1m_1+\Sigma_0m_2>4(r_1+\Sigma_0)n^2,
$$
which implies that $m_1=\mathop{\rm mult}_oZ>4n^2$. Since
$\mathop{\rm deg}Z=4n^2$, we get a contradiction, completing the
proof of Theorem 2.1.


\begin{flushleft}
Department of Mathematical Sciences,\\
The University of Liverpool
\end{flushleft}

\noindent{\it pukh@liverpool.ac.uk}

\end{document}